\documentclass[12pt]{article}
\usepackage{amssymb,amsfonts,amsmath, psfrag,eepic}

\newtheorem{theo}{Theorem}
\newtheorem{defi}[theo]{Definition}

\newtheorem{coro}[theo]{Corollary}

\makeatletter \@addtoreset{equation}{section}
\@addtoreset{theo}{section} \makeatother

\def\qed{\hfill \rule{4pt}{7pt}}
\def\pf{\noindent {\it Proof.} }

\makeatletter
\newfont{\footsc}{cmcsc10 at 8truept}
\newfont{\footbf}{cmbx10 at 8truept}
\newfont{\footrm}{cmr10 at 10truept}

\makeatother 

\begin{document}

\begin{center}

{\Large\bf Riordan Paths and Derangements}

William Y. C. Chen$^1$, Eva Y. P. Deng$^2$ and Laura L. M. Yang$^3$\\[3pt]
\small $^{1,2,3}$Center for Combinatorics, LPMC, Nankai University,\\
\small Tianjin 300071, P. R. China\\[3pt]
\small $^2$Department of Applied Mathematics, Dalian University of Technology, \\
\small Dalian, Liaoning 116024, P.R. China\\[3pt]
\small \texttt{$^1$chen@nankai.edu.cn,  $^2$dengyp@eyou.com,
$^3$yanglm@hotmail.com}

\end{center}

\noindent {\bf Abstract.}  Riordan paths are Motzkin paths without
horizontal steps on the $x$-axis. We establish a correspondence
between Riordan paths and
 $(321,3\bar{1}42)$-avoiding
derangements.  We also present a combinatorial proof of a
recurrence relation for the Riordan numbers in the spirit of the
Foata-Zeilberger proof of a recurrence relation on the Schr\"oder
numbers.

\vskip 8pt

\noindent {\bf Keywords:} Riordan number, Riordan path,
$(321,3\bar{1}42)$-avoiding derangement.

\section{Introduction}

The \emph{Riordan numbers} have many combinatorial
interpretations, see \cite{B} and the On-Line Encyclopedia of
Integer Sequences \cite[A005043]{Sloane}. For example,  the $n$-th
Riordan number $r_n$ equals the number of plane trees with $n$
edges in which no vertex has outdegree one, which are called
\emph{short bushes}. Let $\mathcal{B}_n$ denote the set of short
bushes with $n$ edges (see Figure \ref{table-bush}).  The first
few Riordan numbers are $1, 0, 1, 1, 3, 6, 15, 36, 91, 232.$ In
general, $r_n$ is given by the formula
\begin{equation}
\label{the expression of Riordan number}
 r_n = \frac{1}{n+1}
           \sum_{k=1}^{n-1} {n+1 \choose k}{n-k-1 \choose {k-1}},
\end{equation}
see \cite[A005043]{Sloane}.

The first result of this paper was motivated by the question of
finding a combinatorial interpretation of the Riordan numbers in
terms of permutations with forbidden patterns. In this aspect, we
find that the Riordan numbers are closely related to the Motzkin
numbers. The authors have obtained a combinatorial proof of the
fact that permutations avoiding the patterns $(321,3\bar{1}42)$
are counted by the Motzkin numbers. In this paper, we show that
the Riordan number $r_n$ equals the number of derangements on
$[n]=\{1,2,\ldots,n\}$ that avoid the patterns $(321,3\bar{1}42)$.
Thus the Riordan numbers can be considered as a derangement
analogue of the Motzkin numbers.

The second result of this paper is a combinatorial proof of a
recurrence relation on the Riordan numbers in the spirit of the
Foata-Zeilberger proof of a recurrence on the Schr\"{o}der numbers
\cite{Foata}, see also \cite{Sulanke-1998, Sulanke-1999,
Sulanke-2001}.

\section{Riordan paths}

In this section, we give a brief review of the Riordan numbers and
the Riordan paths. We first give a combinatorial derivation of the
formula \eqref{the expression of Riordan number} by using the
decomposition algorithm obtained in \cite{C}. Let $F_n$ be the set
of labelled plane trees with $n$ edges in which no vertex has
outdegree one. Moreover, let $F_{n,k}$ be the set of trees in
$F_n$ with $k$ internal vertices. Suppose that the set of children
of each internal vertex forms a block.  Using the decomposition
algorithm in \cite{C}, we obtain a bijection between $F_{n,k}$ and
the set of forests with $k$  small plane trees with $n+k$ vertices
such that the roots of the small trees belong to
$\{1,2,\ldots,n+1\}$, and each small tree contains at least two
children. Recall that a small tree is a tree containing only the
root and at least one child.  So $|F_{n,k}|$ can be computed as
follows: we have $\binom{n+1}{k}$ choices for the roots, and the
remaining $n$ different labels are partitioned into $k$ blocks
with each block containing at least two elements. Thus we have
$$
|F_{n,k}|= {n+1 \choose k}{n-k-1 \choose {k-1}}n!,
$$
which implies the formula \eqref{the expression of Riordan number}
because of the relation $F_n = (n+1)! r_n$.

 Recall that a \emph{Motzkin path} of length $n$ is a
lattice path in the plane from $(0,0)$ to $(n,0)$, consisting of
up steps $U=(1,1)$, down steps $D=(1,-1)$, and horizontal steps
$H=(1,0)$, and never going below the $x$-axis \cite{B,Do-S,S}. The
\emph{height} of any step is defined to be the $y$-coordinate of
its starting point. A \emph{$2$-Motzkin path} is a Motzkin path
where the horizontal steps can be of two kinds: straight or wavy.
Motzkin paths are counted by the Motzkin numbers
\cite[A001006]{Sloane} and $2$-Motzkin paths are counted by the
Catalan numbers \cite[A000108]{Sloane}; see, for example,
\cite{De-4, Do-S}.

The Riordan number $r_n$ counts  Motzkin paths of length $n$ with
no horizontal steps of height $0$   \cite[A005043]{Sloane}. This
fact follows from a bijection of Deutsch and Shapiro between plane
trees and $2$-Motzkin paths \cite{De-4}. For any short bush $T$,
let the leftmost  and rightmost edges of a vertex correspond to up
and down steps, respectively, and let the remaining edges
correspond to horizontal steps. Then we obtain a Motzkin path
without horizontal steps on the $x$-axis by traversing $T$ in
preorder.

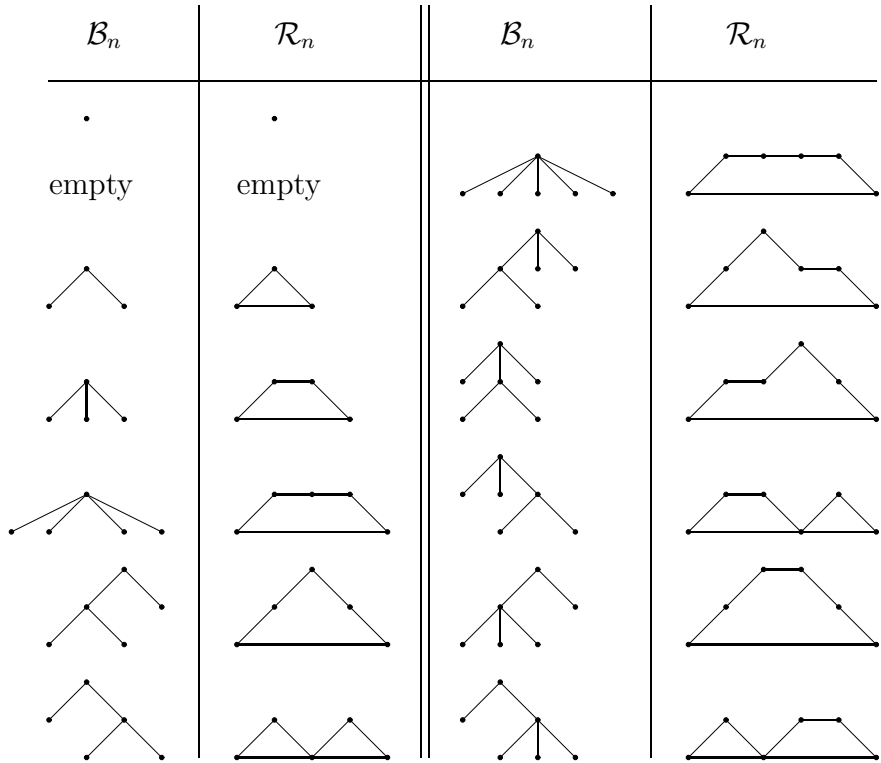
\begin{figure}[h,t]
\begin{center}
\begin{picture}(320,300)
\setlength{\unitlength}{1mm}

\put(15,0){\line(-1,1){5}} \put(5,0){\line(1,1){5}}
\put(10,5){\line(-1,1){5}} \put(0,5){\line(1,1){5}}
\put(5,10){\circle*{0.6}} \put(0,5){\circle*{0.6}}
\put(10,5){\circle*{0.6}} \put(15,0){\circle*{0.6}}
\put(5,0){\circle*{0.6}}

\put(20,0){\line(0,1){100}}

\put(25,0){\line(1,0){20}} \put(25,0){\line(1,1){5}}
\put(30,5){\line(1,-1){5}} \put(35,0){\line(1,1){5}}
\put(40,5){\line(1,-1){5}} \put(25,0){\circle*{0.6}}
\put(30,5){\circle*{0.6}} \put(35,0){\circle*{0.6}}
\put(40,5){\circle*{0.6}} \put(45,0){\circle*{0.6}}

\put(49.5,0){\line(0,1){100}} \put(50.5,0){\line(0,1){100}}

\put(60,0){\line(1,1){5}} \put(60,0){\circle*{0.6}}
\put(65,0){\circle*{0.6}} \put(65,0){\line(0,1){5}}
\put(70,0){\line(-1,1){5}} \put(70,0){\circle*{0.6}}
\put(65,5){\line(-1,1){5}} \put(65,5){\circle*{0.6}}
\put(55,5){\line(1,1){5}} \put(55,5){\circle*{0.6}}
\put(60,10){\circle*{0.6}}

\put(80,0){\line(0,1){100}}

\put(85,0){\line(1,0){25}} \put(85,0){\line(1,1){5}}
\put(85,0){\circle*{0.6}} \put(90,5){\line(1,-1){5}}
\put(90,5){\circle*{0.6}} \put(95,0){\line(1,1){5}}
\put(95,0){\circle*{0.6}} \put(100,5){\line(1,0){5}}
\put(100,5){\circle*{0.6}} \put(105,5){\line(1,-1){5}}
\put(105,5){\circle*{0.6}} \put(110,0){\circle*{0.6}}

\put(0,15){\line(1,1){5}} \put(0,15){\circle*{0.6}}
\put(10,15){\line(-1,1){5}} \put(10,15){\circle*{0.6}}
\put(5,20){\line(1,1){5}} \put(5,20){\circle*{0.6}}
\put(15,20){\line(-1,1){5}} \put(15,20){\circle*{0.6}}
\put(10,25){\circle*{0.6}}

\put(25,15){\line(1,0){20}} \put(25,15){\line(1,1){5}}
\put(25,15){\circle*{0.6}} \put(45,15){\circle*{0.6}}
\put(30,20){\line(1,1){5}} \put(30,20){\circle*{0.6}}
\put(35,25){\line(1,-1){5}} \put(35,25){\circle*{0.6}}
\put(40,20){\line(1,-1){5}} \put(40,20){\circle*{0.6}}

\put(55,15){\line(1,1){5}} \put(55,15){\circle*{0.6}}
\put(60,15){\line(0,1){5}} \put(60,15){\circle*{0.6}}
\put(65,15){\line(-1,1){5}} \put(65,15){\circle*{0.6}}
\put(60,20){\line(1,1){5}} \put(60,20){\circle*{0.6}}
\put(70,20){\line(-1,1){5}} \put(70,20){\circle*{0.6}}
\put(65,25){\circle*{0.6}}

\put(85,15){\line(1,0){25}} \put(110,15){\circle*{0.6}}
\put(85,15){\line(1,1){5}} \put(85,15){\circle*{0.6}}
\put(90,20){\line(1,1){5}} \put(90,20){\circle*{0.6}}
\put(95,25){\line(1,0){5}} \put(95,25){\circle*{0.6}}
\put(100,25){\line(1,-1){5}} \put(100,25){\circle*{0.6}}
\put(105,20){\line(1,-1){5}} \put(105,20){\circle*{0.6}}

\put(-5,30){\line(2,1){10}} \put(-5,30){\circle*{0.6}}
\put(0,30){\line(1,1){5}} \put(0,30){\circle*{0.6}}
\put(10,30){\line(-1,1){5}} \put(10,30){\circle*{0.6}}
\put(15,30){\line(-2,1){10}} \put(15,30){\circle*{0.6}}
\put(5,35){\circle*{0.6}}

\put(25,30){\line(1,0){20}} \put(45,30){\circle*{0.6}}
\put(25,30){\line(1,1){5}} \put(25,30){\circle*{0.6}}
\put(30,35){\line(1,0){5}} \put(30,35){\circle*{0.6}}
\put(35,35){\line(1,0){5}} \put(35,35){\circle*{0.6}}
\put(40,35){\line(1,-1){5}} \put(40,35){\circle*{0.6}}

\put(55,35){\line(1,1){5}} \put(55,35){\circle*{0.6}}
\put(60,35){\line(0,1){5}} \put(60,35){\circle*{0.6}}
\put(65,35){\line(-1,1){5}} \put(65,35){\circle*{0.6}}
\put(60,30){\line(1,1){5}} \put(60,30){\circle*{0.6}}
\put(70,30){\line(-1,1){5}} \put(70,30){\circle*{0.6}}
\put(60,40){\circle*{0.6}}

\put(85,30){\line(1,0){25}} \put(110,30){\circle*{0.6}}
\put(85,30){\line(1,1){5}} \put(85,30){\circle*{0.6}}
\put(90,35){\line(1,0){5}} \put(90,35){\circle*{0.6}}
\put(95,35){\line(1,-1){5}} \put(95,35){\circle*{0.6}}
\put(100,30){\line(1,1){5}} \put(100,30){\circle*{0.6}}
\put(105,35){\line(1,-1){5}} \put(105,35){\circle*{0.6}}

\put(0,45){\line(1,1){5}} \put(0,45){\circle*{0.6}}
\put(5,45){\line(0,1){5}} \put(5,45){\circle*{0.6}}
\put(10,45){\line(-1,1){5}} \put(10,45){\circle*{0.6}}
\put(5,50){\circle*{0.6}}

\put(25,45){\line(1,0){15}} \put(40,45){\circle*{0.6}}
\put(25,45){\line(1,1){5}} \put(25,45){\circle*{0.6}}
\put(30,50){\line(1,0){5}} \put(30,50){\circle*{0.6}}
\put(35,50){\line(1,-1){5}} \put(35,50){\circle*{0.6}}

\put(55,45){\line(1,1){5}} \put(55,45){\circle*{0.6}}
\put(65,45){\line(-1,1){5}} \put(65,45){\circle*{0.6}}
\put(60,50){\line(0,1){5}} \put(60,50){\circle*{0.6}}
\put(65,50){\line(-1,1){5}} \put(65,50){\circle*{0.6}}
\put(55,50){\line(1,1){5}} \put(55,50){\circle*{0.6}}
\put(60,55){\circle*{0.6}}

\put(85,45){\line(1,0){25}} \put(110,45){\circle*{0.6}}
\put(85,45){\line(1,1){5}} \put(85,45){\circle*{0.6}}
\put(90,50){\line(1,0){5}} \put(90,50){\circle*{0.6}}
\put(95,50){\line(1,1){5}} \put(95,50){\circle*{0.6}}
\put(100,55){\line(1,-1){5}} \put(100,55){\circle*{0.6}}
\put(105,50){\line(1,-1){5}} \put(105,50){\circle*{0.6}}

\put(0,60){\line(1,1){5}} \put(0,60){\circle*{0.6}}
\put(10,60){\line(-1,1){5}} \put(10,60){\circle*{0.6}}
\put(5,65){\circle*{0.6}}

\put(25,60){\line(1,0){10}} \put(35,60){\circle*{0.6}}
\put(25,60){\line(1,1){5}} \put(25,60){\circle*{0.6}}
\put(30,65){\line(1,-1){5}} \put(30,65){\circle*{0.6}}

\put(55,60){\line(1,1){5}} \put(55,60){\circle*{0.6}}
\put(65,60){\line(-1,1){5}} \put(65,60){\circle*{0.6}}
\put(60,65){\line(1,1){5}} \put(60,65){\circle*{0.6}}
\put(65,65){\line(0,1){5}} \put(65,65){\circle*{0.6}}
\put(70,65){\line(-1,1){5}} \put(70,65){\circle*{0.6}}
\put(65,70){\circle*{0.6}}

\put(85,60){\line(1,0){25}} \put(110,60){\circle*{0.6}}
\put(85,60){\line(1,1){5}} \put(85,60){\circle*{0.6}}
\put(90,65){\line(1,1){5}} \put(90,65){\circle*{0.6}}
\put(95,70){\line(1,-1){5}} \put(95,70){\circle*{0.6}}
\put(100,65){\line(1,0){5}} \put(100,65){\circle*{0.6}}
\put(105,65){\line(1,-1){5}} \put(105,65){\circle*{0.6}}

\put(0,75){empty}

\put(25,75){empty}

\put(5,85){\circle*{0.6}}

\put(30,85){\circle*{0.6}}

\put(55,75){\line(2,1){10}} \put(55,75){\circle*{0.6}}
\put(60,75){\line(1,1){5}} \put(60,75){\circle*{0.6}}
\put(65,75){\line(0,1){5}} \put(65,75){\circle*{0.6}}
\put(70,75){\line(-1,1){5}} \put(70,75){\circle*{0.6}}
\put(75,75){\line(-2,1){10}} \put(75,75){\circle*{0.6}}
\put(65,80){\circle*{0.6}}

\put(85,75){\line(1,0){25}} \put(110,75){\circle*{0.6}}
\put(85,75){\line(1,1){5}} \put(85,75){\circle*{0.6}}
\put(90,80){\line(1,0){5}} \put(90,80){\circle*{0.6}}
\put(95,80){\line(1,0){5}} \put(95,80){\circle*{0.6}}
\put(100,80){\line(1,0){5}} \put(100,80){\circle*{0.6}}
\put(105,80){\line(1,-1){5}} \put(105,80){\circle*{0.6}}

\put(0,90){\line(1,0){110}}

\put(5,95){\bf $\mathcal{B}_n$}

\put(30,95){\bf $\mathcal{R}_n$}

\put(60,95){\bf $\mathcal{B}_n$}

\put(90,95){\bf $\mathcal{R}_n$}
\end{picture}
\end{center}
\caption{Short bushes and Riordan paths} \label{table-bush}
\end{figure}

A Motzkin path of length $n$ without horizontal steps on the
$x$-axis will be called a {\it Riordan path} of length $n$, and
let $\mathcal{R}_n$ be the set of Riordan paths of length $n$.
Figure \ref{table-bush} is an illustration of the  correspondence
between short bushes and Riordan paths.

The Riordan numbers $r_n$ are related to the Catalan numbers
$c_n=\frac{1}{n+1}\binom{2n}{n}$ by the relation
     \begin{equation}
        \label{re-cr}
         c_n=\sum_{k=0}^{n} {n \choose k} r_k,
     \end{equation}
which leads to the following formula:
\begin{equation}
        \label{re-rc}
        r_{n}=\sum_{k=0}^n(-1)^{n-k}{n \choose k}c_k.
\end{equation}
The above  formula (\ref{re-rc}) has been derived by Bernhart
\cite{B} using a difference operator. Here we present a
combinatorial interpretation of (\ref{re-cr}).

\noindent {\it Combinatorial Proof of (\ref{re-cr}).}
 Let $P=p_1 p_2 \cdots p_{2n}$ be a Dyck path of length $2n$.
 We divide the path $P$ into $n$
segments $Q_1Q_2\cdots Q_n$ such that $Q_{i}=p_{2i-1}p_{2i}$. For
each $Q_i$, there are four possible combinations: $UU$, $UD$, $DU$
and $DD$. If we use the four kinds of steps of a $2$-Motzkin path
to encode $UU$, $UD$, $DU$ and $DD$, that is, $UU$ is represented
by an up step, $UD$ is represented by a wavy horizontal step, $DU$
is represented by a straight horizontal step, and $DD$ is
represented by a down step. Then we get a $2$-Motzkin path $M$
without  straight horizontal steps on the $x$-axis. Suppose $M$
contains $n-k$ wavy horizontal steps. Note that if we remove all
the wavy horizontal steps, we are led to a Riordan path of length
$k$. Conversely, given a Riordan path of length $k$, we can
reconstruct ${n \choose k}$ $2$-Motzkin paths without straight
horizontal steps on the $x$-axis by inserting $n-k$ wavy
horizontal steps. \qed

The above proof implies the following interpretation of the
Catalan number $c_n={1 \over n+1} {2n \choose n}$.

\begin{coro}
The number of  $2$-Motzkin paths of length $n$ without straight
horizontal steps on the $x$-axis equals the Catalan number $c_n$.
\end{coro}

\section{Riordan Paths and Derangements}

In this section, we give a correspondence between Riordan paths
and derangements with forbidden patterns $(321,3\bar{1}42)$. This
is motivated by the recent work of the authors \cite{C-D-Y} on the
bijection  $\phi$ between Motzkin paths of length $n$ and
$S_n(321, 3\bar{1}42)$, where $S_n$ denotes the set of
permutations on $[n]$, and $S_n(321, 3\bar{1}42)$ denote the set
of permutations avoiding the patterns $(321,3\bar{1}42)$. We say
that a permutation $\pi=\pi_1\pi_2\cdots \pi_n$ avoids the pattern
$321$ if it does not contain any subsequence $\pi_i\pi_j\pi_k$
such that $\pi_i> \pi_j>\pi_k$ for $1 \leq i<j<k\leq n$. Moreover,
we say that $\pi$ avoids the pattern $3\bar{1}42$ if any
subsequence $\pi_i\pi_j\pi_k$ $(i<j<k)$ of pattern $231$, namely,
$\pi_j> \pi_i>\pi_k$, can be extended to a subsequence of pattern
$3142$, in other words, there exists $i<m<j$ such that $\pi_j>
\pi_i>\pi_k>\pi_m$.

It was shown by Gire \cite{Gire} that $|S_n(321, 3\bar{1}42)|$
equals the Motzkin number $m_n$ (see \cite[A001006]{Sloane}).
Authors \cite{C-D-Y} established a correspondence between Motzkin
paths of length $n$ and reduced decompositions of permutations in
$S_n(321, 3\bar{1}42)$. In order to make a connection between
Riordan paths and permutations with forbidden patterns, we led to
the consideration of further restrictions on $S_n(321,
3\bar{1}42)$ so that we may get a subset of permutations $S_n(321,
3\bar{1}42)$ that are in one-to-one correspondence with Riordan
paths of length $n$ with $m$ horizontal steps on the $x$-axis.

We now recall the definition of $\phi$ which is given in terms of
reduced decompositions of permutations in $S_n$.

\begin{defi}
For any $1\leq i\leq n-1$, define the map $s_i$: ${{S}}_n
\rightarrow {{S}}_n$, such that $s_i$ acts on a permutation by
interchanging the elements in positions $i$ and $i+1$. We call
$s_i$ the simple transposition, and write the action of $s_i$ on
the right of the permutation, denoted by $\pi s_i$. Therefore we
have $\pi (s_i s_j)=(\pi s_i)s_j$.
\end{defi}

The canonical reduced decomposition of $\pi \in {S}_n$ has the
following form:
\begin{equation}\label{dck}
\pi=(1\ 2\ \cdots\ n)\sigma = (1 \ 2\ \cdots \ n)
\sigma_1\sigma_2\cdots\sigma_k,
\end{equation}
 where
\begin{equation*}
 \sigma_i= s_{h_i} s_{h_i-1} \cdots s_{t_i},  \quad
  h_i\geq t_i \quad (1\leq i \leq k)\quad \mbox{and}
  \end{equation*}
  \[
   1\leq h_1< h_2< \cdots < h_k\leq n-1.
   \]
 We call $h_i$ the
{\it head} and $t_i$ the {\it tail} of $\sigma_i$. For short, we
say that $\pi$ has the canonical reduced decomposition
$\sigma_1\sigma_2\cdots \sigma_k$.

For example,  $\pi=315264$ has the canonical reduced decomposition
$(s_2s_1)(s_4s_3)(s_5)$. It is shown in \cite{C-D-Y} that
permutations in $S_n(321, 3\bar{1}42)$ can be characterized by
their reduced decompositions.

\begin{theo}\label{peravoid}
Let $\pi$ be a permutation in $S_n$ with the reduced decomposition
as given in (\ref{dck}). Then $\pi\in S_n(321, 3\bar{1}42)$ if and
only if
\begin{equation}\label{condition-decomposition}
  t_{i+1} \geq t_{i}+2,\ \ \quad 1\leq i\leq k-1.
\end{equation}
\end{theo}

We now give a brief description of the bijection $\phi$ between
Motzkin paths of length $n$ and $S_n(321,3\bar{1}42)$ by the
\emph{strip decomposition} of Motzkin paths \cite{C-D-Y}. This
bijection involves a labelling of the cells in the region of a
Motzkin path. The region of a Motzkin path is meant to be the area
surrounded by the path and the $x$-axis. Furthermore, the region
of a Motzkin path is subdivided into cells which are either unit
squares or triangles with unit bottom sides. A triangular cell
contains either an up step or a down step. We will not label
triangular cells containing up steps. The other types of cells,
either square or triangular, have bottom sides, say, with  points
$(i,j)$ and $(i+1, j)$, we will label these cells with $s_{i+j}$
or simply $i+j$. We call this labelling the $(x+y)$-labelling.

We now define the strip decomposition of a Motzkin path. Suppose
$P_{n,k}$ is a Motzkin path of length $n$ that contains $k$ up
steps. If $k=0$, then the strip decomposition of $P_{n,0}$ is
simply the empty set.  For any $P_{n,k} \in {M}_n$, let
$A\rightarrow B$ be the last up step and $ E \rightarrow F$ the
last down step on $P_{n,k}$. Then we define the strip of $P_{n,k}$
as the path from $B$ to $F$ along the path $P_{n,k}$. Now we move
the points from $B$ to $E$ one layer lower, namely, subtract the
$y$-coordinate by $1$, and denote the adjusted points by $B'$,
$\ldots$, $E'$. We now form a new Motzkin path by using the path
$P_{n,k}$ up to the point $A$, then joining the point $A$ to $B'$
and following the adjusted segment until we reach the point $E'$,
then continuing with the points on the $x$-axis to reach the
destination $(n, 0)$. Denote this Motzkin path by $P_{n,k-1}$,
which may end with some horizontal steps.

From the strip of $P_{n,k}$, we may define the value $h_k$ as the
label of the cell containing the step $E\rightarrow F$. Clearly,
we have $h_k\leq n-1$. The value $t_k$ is defined as the label of
the cell containing the step starting from the point $B$.

\begin{figure}[h,t]
\begin{center}
\begin{picture}(350,80)
\setlength{\unitlength}{1.5mm} \put(-1,0){\line(1,0){86}}
\put(0,-1){\line(0,1){21}}

\multiput(0,5)(1,0){85}{\line(1,0){0.7}}
\multiput(0,10)(1,0){85}{\line(1,0){0.7}}
\multiput(0,15)(1,0){85}{\line(1,0){0.7}}
\multiput(0,20)(1,0){85}{\line(1,0){0.7}}

\multiput(5,0)(0,1){20}{\line(0,1){0.7}}
\multiput(10,0)(0,1){20}{\line(0,1){0.7}}
\multiput(15,0)(0,1){20}{\line(0,1){0.7}}
\multiput(20,0)(0,1){20}{\line(0,1){0.7}}
\multiput(25,0)(0,1){20}{\line(0,1){0.7}}
\multiput(30,0)(0,1){20}{\line(0,1){0.7}}
\multiput(35,0)(0,1){20}{\line(0,1){0.7}}
\multiput(40,0)(0,1){20}{\line(0,1){0.7}}
\multiput(45,0)(0,1){20}{\line(0,1){0.7}}
\multiput(50,0)(0,1){20}{\line(0,1){0.7}}
\multiput(55,0)(0,1){20}{\line(0,1){0.7}}
\multiput(60,0)(0,1){20}{\line(0,1){0.7}}
\multiput(65,0)(0,1){20}{\line(0,1){0.7}}
\multiput(70,0)(0,1){20}{\line(0,1){0.7}}
\multiput(75,0)(0,1){20}{\line(0,1){0.7}}
\multiput(80,0)(0,1){20}{\line(0,1){0.7}}
\multiput(85,0)(0,1){20}{\line(0,1){0.7}}

\thicklines \put(0,0){\line(1,1){5}} \put(5,5){\line(1,0){5}}
\put(10,5){\line(1,0){5}} \put(15,5){\line(1,-1){5}}
\put(20,0){\line(1,1){5}} \put(25,5){\line(1,1){5}}
\put(30,10){\line(1,0){10}} \put(40,10){\line(1,-1){5}}
\put(45,5){\line(1,0){5}} \put(50,5){\line(1,1){10}}
\put(60,15){\line(1,-1){10}} \put(70,5){\line(1,0){10}}
\put(80,5){\line(1,-1){5}}

\put(7,1){1} \put(12,1){2} \put(17,1){3} \put(27,1){5}
\put(32,1){6} \put(37,1){7} \put(42,1){8} \put(47,1){9}
\put(51,1){10} \put(56,1){11} \put(61,1){12} \put(66,1){13}
\put(71,1){14} \put(76,1){15} \put(81,1){16} \put(32,6){7}
\put(37,6){8} \put(42,6){9} \put(56,6){12} \put(61,6){13}
\put(66,6){14} \put(61,11){14} \put(66,6){14} \thicklines

\put(55,10){\line(1,0){5}} \put(60,10){\line(1,-1){10}} \thinlines
\put(70,0){\line(1,0){10}} \thicklines \put(50,5){\line(1,0){10}}
\put(60,5){\line(1,-1){5}} \put(25,5){\line(1,0){15}}
\put(40,5){\line(1,-1){5}}

\put(-3,-3){(0,0)}
\end{picture}
\end{center}
\caption{The $(x+y)$-labeling and strip decomposition}
\label{labels}
\end{figure}
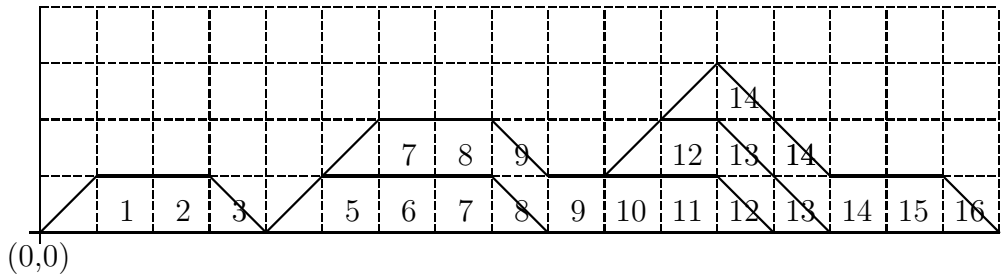

Iterating the above procedure, we get a set of parameters $\{(h_i,
t_i) | 1 \leq i \leq k\}$ satisfying the condition
\eqref{condition-decomposition}. For each step in the above
procedure, we obtain a product of transpositions $\sigma_i=
s_{h_i} s_{h_i-1} \cdots s_{t_i}$. Finally we get the
corresponding canonical reduced decomposition
$\sigma=\sigma_1\sigma_2\cdots \sigma_k$ and the corresponding
permutation $\pi=(1~2~\cdots~n)\sigma$, see Figure \ref{labels}.
We then obtain the following property of the bijection $\phi$.

\begin{theo}
The bijection $\phi$ is a correspondence between Motzkin paths of length $n$ with $m$
horizontal steps on the $x$-axis and permutations in $S_n(321, 3\bar{1}42)$ that have $m$
fixed points.
\end{theo}
\pf For any Motzkin path $P$ of length $n$ with $m$ horizontal
steps on the $x$-axis, label its steps with $0,1,2,\ldots, n-1$
from left to right. Suppose that the $m$ horizontal steps on the
$x$-axis are labelled by $x_1, x_2, \ldots, x_m$, where $0\leq
x_1<x_2<\ldots <x_m\leq n-1$. By the strip decomposition and the
$(x+y)$-labelling, $s_{x_1}, s_{x_2},\ldots, s_{x_m}$ do not occur
in its corresponding canonical reduced decomposition with respect
to the bijection $\phi$. Note that a horizontal step on the
$x$-axis is followed by an up step or a horizontal step on the
$x$-axis (except that it is the last step). Thus $x_1+1, x_2+1,
\ldots, x_m+1$  are  fixed points of the corresponding permutation
in $S_n(321, 3\bar{1}42)$ by applying Theorem \ref{peravoid}. \qed

\begin{coro}
For any Motzkin path $P$ of length $n$, let $\pi \in
S_n(321,3\bar{1}42)$ be its corresponding permutation with respect
to the bijection $\phi$. Suppose that $\pi$ has the canonical
reduced decomposition of the form \eqref{dck}, then
\begin{itemize}
    \item [{\rm 1)}] $t_1-1$ is the number of initial horizontal steps
    on the $x$-axis at the beginning of the Motzkin path $P$;
    \item [{\rm 2)}] $n-1-h_k$ is the number of final horizontal steps on the
    $x$-axis at the end of the Motzkin path $P$;
    \item [{\rm 3)}] $\sum_{i}(t_{i+1}-h_i-2)$ equals the number of
     horizontal steps  of the Motzkin
    path $P$ on the $x$-axis that are neither initial nor final
    steps, where summation is over all $i$ such that
    $h_i+1<t_{i+1}$.
\end{itemize}
\end{coro}

Recall that a permutation $\pi=\pi_1\pi_2\cdots \pi_n$ is said to
be a {\it derangement} if $\pi$ does not have any fixed points,
that is, $\pi_i\neq i$ for all $i \in [n]$. Let $D_n(321,
3\bar{1}42)$ denote $(321, 3\bar{1}42)$-avoiding derangements in
${S}_n$. Then we have the following correspondence.
\begin{coro}
The bijection $\phi$ is a correspondence between Riordan paths of
length $n$ and $D_n(321, 3\bar{1}42)$.
\end{coro}

For example, for the Riordan path in Figure \ref{labels}, we have
\[P_{17,5}=UHHDUUHHDHUUDDHHD.\]
From the strip decomposition, we get the parameter set
\[\{(3, 1), \; (8,5), \; (12, 7), \; (13, 12), \; (16,14)\}\]
The canonical reduced decomposition is given below:
\begin{equation}
\label{eperm}
(s_3s_2s_1)(s_8s_7s_6s_5)(s_{12}s_{11}s_{10}s_9s_8s_7)
(s_{13}s_{12})(s_{16}s_{15}s_{14}).
\end{equation}
The corresponding permutation is
\[ 4 \ \ 1 \ \ 2 \ \ 3\ \ 9\ \ 5\ \ 13\ \ 6\ \ 7\ \
8\ \ 10\  \ 14\ \ 11\ \ 17\ \ 12\ \ 15\  \ 16.\]

\begin{coro}
Let $P$ be a Riordan path of length $n$. Then the area of $P$
minus the sum of heights of  the up steps is equal to the
inversion number of the permutation $\phi(P)\in
D_n(321,3\bar{1}42)$.
\end{coro}

\begin{coro}\label{main} Let $\sigma=\sigma_1\cdots\sigma_k$ be the canonical
reduced decomposition of $\pi\in {S}_n$, where $\sigma_i = s_{h_i}
s_{h_i-1} \cdots s_{t_i}$ for $1\leq i\leq k$. Then $\pi \in
D_n(321, 3\bar{1}42)$ if and only if $t_1=1$, $h_k=n-1$ and
$$
  h_i+2 \geq t_{i+1}  \geq t_{i} + 2,\ \ \quad 1\leq i \leq k-1.
$$
\end{coro}

\section{A Recurrence Relation }

In this section, we give a combinatorial proof of the following
recurrence relation on the Riordan numbers:

\begin{theo} For $n\geq 2$, we have
\begin{equation}\label{re-rr}
(n+1) r_n = (n-1)(2r_{n-1}+3r_{n-2}),
\end{equation}
with initial values $r_0=1$, $r_1=0$ and $r_2=1$.
\end{theo}
\pf  We proceed to establish the following bijection:
\begin{equation}
\label{map-e}
 \psi\colon \ \ [3(n-1)]\times\mathcal{R}_{n-2}\ \bigcup\;\; [2(n-1)] \times \mathcal{R}_{n-1}
        \;\; \Longrightarrow \;\;
          [n+1] \times \mathcal{R}_n,
\end{equation}
which yields the identity (\ref{re-rr}).

We begin with an interpretation of  $[3(n-1)]\times
\mathcal{R}_{n-2}$ as the multi-set of Riordan paths of length
$n-2$ in which exactly one step is labelled one of the labels $a$,
$b$, and $c$, plus three copies of the set of Riordan paths of
length $n-2$ without labels. Similarly,  $[2(n-1)] \times
\mathcal{R}_{n-1}$ can be represented by the set of labelled
Riordan paths of length $n-1$ in which exactly one step is
labelled either by $1$ or $2$. The set $[n+1]\times \mathcal{R}_n$
can be represented by the set of Riordan paths of length $n$ for
which at most one step is labelled by the symbol $*$.

For example, since $\mathcal{R}_4=\{UUDD, UDUD, UHHD\}$,
$[5]\times \mathcal{R}_4$ consists of the following labelled
paths:
\begin{equation*}
   \begin{array}{rllll}
        UUDD & U^* UDD & UU^* DD & UUD^* D & UUDD^* \\
        UDUD & U^* DUD & UD^* UD & UDU^* D & UDUD^* \\
        UHHD & U^* HHD & UH^* HD & UHH^* D & UHHD^* .
          \end{array}
\end{equation*}

We  now  give a construction of the map $\psi$.

(1) For the three copies of the paths in $\mathcal{R}_{n-2}$
without labels, we  respectively add $UD$, $U^*D$ and $UD^*$  to
the beginning of the paths. In this way, we obtain all the paths
beginning with $UD$ in $[n+1]\times \mathcal{R}_n$. For example,
for $n=4$, the three copies of $UD$ are mapped to $UDUD$, $U^*DUD$
and $UD^*UD$, respectively.

(2) For the paths having a step $p_i$ of height $k$ labelled by
 $a$ in $\mathcal{R}_{n-2}$: If $k=0$, namely, $p_i=U$,
we add an up step to the beginning of the path and insert a down
step following the corresponding down step of $p_i$, namely, the
first down step after $p_i$ that touches the $x$-axis. This gives
all the Riordan paths of length $n$ without labels such that there
is no horizontal steps of height $1$ before the path returns to
the $x$-axis. Otherwise, let $p_j$ be the last up step of height
$k-1$ before the step $p_i$, then we add an up step after $p_j$
and a down step before $p_i$ and label $p_j$ with $*$. Hence we
have all the Riordan paths of length $n$ which contain the
consecutive steps $U^*U$. For example, $U^aD$ and $UD^a$ are
mapped to $UUDD$ and $U^*UDD$, respectively.

(3) For the paths having a step $p_i$ labelled by $b$ (or $c$) in
$\mathcal{R}_{n-2}$, we add $U^*D$ (or $UD^*$) after $p_i$. In
this way, we get all Riordan paths of length $n$ containing the
consecutive steps $U^*D$ (or $UD^*$) which are not at the
beginning of the Riordan paths. For example, $U^bD$ and $UD^b$ (or
$U^cD$ and $UD^c$) are mapped to $UU^*DD$ and $UDU^*D$ (or
$UUD^*D$ and $UDUD^*$), respectively.

(4) For the paths having a step $p_i$ of height $k$ labelled by
$1$ in $\mathcal{R}_{n-1}$: If $p_i=D$ and $k=1$, then we change
the corresponding up step (that is, the nearest up step before
$p_i$ that touches the $x$-axis) to an $H$ step, and add an up
step to the beginning of the path. So we obtain all the Riordan
paths of length $n$ without labels such that there is at least one
horizontal step of height $1$ before the path returns to the
$x$-axis. Otherwise, we add a horizontal step after $p_i$, and
label the new horizontal step with $*$. This yields all the
Riordan paths of length $n$ containing $H^*$. For example,
$U^1HD$, $UH^1D$ and $UHD^1$ are mapped to $UH^*HD$, $UHH^*D$ and
$UHHD$, respectively.

(5) For the paths having a step $p_i$ labelled by $2$ in
$\mathcal{R}_{n-1}$: If $p_i$ is an up step (or a down step), then
we label $p_i$ with $*$ and add a horizontal step $H$ after $p_i$
(before $p_i$). Thus we obtain all the Riordan paths of length $n$
containing the consecutive steps $U^*H$ (or  $HD^*$). If $p_i=H$,
then its height is nonzero. In this case,  so we may assume that
$p_j$ is the first down step after $p_i$. Then we replace $p_i$ by
$U$, and add a down step before $p_j$ and label $p_j$ with $*$. So
we obtain all the Riordan paths of length $n$ containing
consecutive steps $DD^*$. For example, $U^2HD$, $UH^2D$ and
$UHD^2$ are mapped to $U^*HHD$, $UUDD^*$, $UHHD^*$, respectively.

In summary, we obtain all the Riordan paths in $[n+1]\times
\mathcal{R}_n$. It can be seen that the above procedure is
reversible. Hence $\psi$ is a bijection. \qed

Note that the relation \eqref{re-rr} is derived from the
generating function by Bernhart \cite{B}. Our proof is in the
spirit  of the Foata-Zeilberger proof of a recurrence relation on
the Schr\"oder numbers \cite{Foata}, and Sulanke's proofs of the
recurrences for Schr\"oder paths, parallelogram polyominoes and
Motzkin paths \cite{Sulanke-1998, Sulanke-1999, Sulanke-2001}.

\vspace{0.5cm} \noindent{\bf Acknowledgments.} We are grateful to
the referees for valuable comments.  This work was supported by
the 973 Project on Mathematical Mechanization, and the National
Science Foundation of China, the Ministry of Education, and the
Ministry of Science and Technology of China.

\end{document}